\documentclass[12pt]{amsart}
\usepackage{amssymb,latexsym}
\usepackage{epsfig}
\usepackage{url}
\usepackage{enumerate}

\newtheorem*{theorem*}{Theorem}

\newcommand{\Rbb}{\mathbb{R}}

\newcommand{\calh}{\mathcal{H}}

\newcommand{\pxone}{\frac{\partial}{\partial x_{1}}}
\newcommand{\pyone}{\frac{\partial}{\partial y_{1}}}
\newcommand{\pxtwo}{\frac{\partial}{\partial x_{2}}}
\newcommand{\pytwo}{\frac{\partial}{\partial y_{2}}}

\newcommand{\mgap}{\vspace{24pt}}
\newcommand{\sgap}{\vspace{12pt}}

\begin{document}

\title{The Length Spectrum of the sub-Riemannian Three-Sphere}
\author{David Klapheck and Michael VanValkenburgh}
\address{Department of Mathematics and Statistics, California State University Sacramento,
Sacramento, 95819, USA}
\email{dtk22@csus.edu, mjv@csus.edu}
\date{\today}

\maketitle

\renewcommand*{\thefootnote}{\fnsymbol{footnote}}

\begin{abstract}
	We determine the lengths of all closed sub-Riemannian geodesics on the three-sphere $S^{3}$. Our methods are elementary and allow us to avoid using explicit formulas for the sub-Riemannian geodesics.
\end{abstract}

\section{Introduction}

In the case of a compact Riemannian manifold $(M,g)$ there is a relationship between closed geodesics, representing paths of free classical particles in periodic motion, and eigenfunctions of the Laplacian $\Delta$, representing periodic free quantum ``waves'' (up to a phase factor). For this reason, the set of lengths of closed geodesics is called the \emph{length spectrum}, in analogy to the spectrum of the Laplacian. There are in fact precise formulas relating lengths to eigenvalues; see for example the announcement of Guillemin and Weinstein for a readable discussion with references \cite{R:GuillWein}.

\mgap

So far there is no such formula relating lengths and eigenvalues in the case of a compact sub-Riemannian (sR) manifold. We recall that a sR manifold is a manifold with a specified linear subbundle $\calh$ (the ``horizontal bundle'') of its tangent bundle, along with a Riemannian metric on $\calh$. Distances between points are then measured using curves that are constrained to have tangent vectors in $\calh$ (``horizontal curves''). In fact, when $\calh$ is the span of a set of bracket-generating vector fields, then the Chow-Rashevskii theorem says that any two points are connected by a horizontal curve, a result that even experts find surprising (\cite{R:Burago}, p.178); thus given any two points there is a shortest horizontal curve connecting them; it is called a \emph{sR geodesic}.

\mgap

sR geometry is of practical interest; for example, the problem of parallel parking a car, or, even worse, a car with a trailer, is a problem in sR geometry \cite{R:Burago},\cite{R:Nelson}. And there are further surprises from the purely mathematical point of view, one being Montgomery's proof of existence of singular sR geodesics (singular in the sense that they do not satisfy the geodesic equations (Hamilton's equations)) \cite{R:Montg94a},\cite{R:Montg}. This and other relatively recent results in sR geometry then inspire renewed interest in the subLaplacian: the operator naturally associated with the given (sub-)Riemannian metric on $\calh$.

\mgap

In this paper, with the goal of understanding a single example, we compute the sR length spectrum of the three-dimensional sphere $S^{3}$ with its standard sR structure; this is to be compared with the spectrum of the subLaplacian on $S^{3}$, known by Taylor \cite{R:TaylorNonc} and generalized to other connected, semisimple Lie groups by Domokos \cite{R:DomokosPeterWeyl}. We expect that a general theory relating the sR length spectrum to the spectrum of the subLaplacian would be amenable to the tools of microlocal analysis as in the Riemannian setting; the work of Colin de Verdi\`{e}re, et al., gives hope that this will be accomplished \cite{R:CdVproc}. 

\mgap

We focus on $S^{3}$ with its standard sR structure because it is perhaps the simplest compact manifold with a sR structure, and there are no singular sR geodesics on $S^{3}$; that is, all sR geodesics arise as projections of solutions of Hamilton's equations \cite{R:Montg}. Moreover, we wish to compare the sR setting to the Riemannian setting, in which the spheres $S^{n}$ are of fundamental importance, as examples of manifolds all of whose geodesics are closed and have the same length $T$; in general this is equivalent to most of the spectrum of $\sqrt{-\Delta}$ being concentrated near an arithmetic progression $\frac{2\pi}{T}k+\beta$, $k=1,2,\ldots$, for some constant $\beta$ \cite{R:DuistGuill}. As we will see, in the case of $S^{3}$ not all sR geodesics are closed, and not all have the same length: 

\sgap

\begin{theorem*}
	The set of lengths of the closed sR geodesics on $S^{3}$ is
	$$\{2\pi\sqrt{n};\,\,n\in\mathbb{N}\}.$$
\end{theorem*}

\mgap

Others have studied the sR geodesics on $S^{3}$ \cite{R:CaChMa}, \cite{R:ChMaVa}, \cite{R:HurtadoRosales} (see also the survey article \cite{R:DAngTyson}), but we compute their lengths and differ from the previous work in that we consistently use Hopf coordinates on $S^{3}$ and avoid using explicit formulas for the sR geodesics; we believe it clarifies the presentation to \emph{not} use explicit formulas.

\mgap

We introduce the sR structure and geodesic equations in Section~\ref{S:Hopf} using Hopf coordinates, and in Section~\ref{S:Categ} we categorize the qualitatively different types of sR geodesics. In Section~\ref{S:closure} we determine which sR geodesics are closed, and in Section~\ref{S:lengths} we compute their lengths, resulting in the theorem above. Finally, in Section~\ref{S:spectrum} we compare the sR length spectrum to the previously-known spectrum of the subLaplacian.

\mgap

This work was supported by a SURE (Summer Undergraduate Research Experience) Award at California State University, Sacramento.

\mgap

\noindent\textbf{Remark.} During peer review, it was pointed out that the above result is contained in the paper \cite{R:ChMaVaHopf} (see their Theorem~2). However, our proof is entirely new and has the advantage of being elementary after the introduction of Hamilton's equations (\ref{E:HamEqHopf}) in our chosen coordinate system.

\newpage

\section{$S^{3}$ in Euclidean and Hopf Coordinates}\label{S:Hopf}

First we consider $S^{3}$ as a subset of $\Rbb^{4}$: 
$$S^{3}=\{(x_{1},y_{1},x_{2},y_{2})\in\Rbb^{4};\,\,x_{1}^{2}+y_{1}^{2}+x_{2}^{2}+y_{2}^{2}=1\}.$$
On $S^{3}$ we have the orthonormal vector fields
\begin{align*}
	V&:=-y_{1}\pxone+x_{1}\pyone-y_{2}\pxtwo+x_{2}\pytwo\\
	E_{1}&:=-x_{2}\pxone+y_{2}\pyone+x_{1}\pxtwo-y_{1}\pytwo\\
	E_{2}&:=-y_{2}\pxone-x_{2}\pyone+y_{1}\pxtwo+x_{1}\pytwo
\end{align*}
which satisfy the Lie bracket relations
\begin{align*}
	[V,E_{1}]=-2E_{2},\qquad [E_{2},V]=-2E_{1},\qquad [E_{1},E_{2}]=-2V.
\end{align*}
Thus $\mathcal{H}(S^{3})=\text{Span}\{E_{1},E_{2}\}$ is a bracket-generating tangent subbundle, and by the Chow-Rashevskii Theorem any two points on $S^{3}$ are connected by a sR geodesic.

\mgap

The orbits of the flow generated by $V$ are the circles of the Hopf fibration \cite{R:CdS}, so we find it convenient to use Hopf coordinates\footnote{``$3$-sphere.'' Wikipedia, The Free Encyclopedia. Wikimedia Foundation, Inc. June 16, 2015. Web. July 1, 2015.} on $S^{3}$:
\begin{equation*}
\begin{aligned}
	x_{1}&=\cos\theta_{1}\sin\theta_{0}\\
	y_{1}&=\sin\theta_{1}\sin\theta_{0}\\
	x_{2}&=\cos\theta_{2}\cos\theta_{0}\\
	y_{2}&=\sin\theta_{2}\cos\theta_{0},
\end{aligned}
\end{equation*}
for $0<\theta_{0}<\pi/2$ and $0<\theta_{j}<2\pi$, $j=1,2$. We picture the $(\theta_{0},\theta_{1},\theta_{2})$-space as ``the Hopf cube'' $(0,\pi/2)\times(0,2\pi)\times(0,2\pi)$. When we have occasion to exit the Hopf cube, we simply return to the definition of Hopf coordinates to make the correct interpretation:

\sgap

\begin{enumerate}[(i)]
	\item For the $\theta_{1}$ and $\theta_{2}$ coordinates the values $0$ and $2\pi$ are identified.
	\item When a point crosses the $\theta_{0}=0$ plane we have that $\theta_{0}$ changes direction (``bounces'') and $(\theta_{1},\theta_{2})$ is identified with $(\theta_{1}+\pi,\theta_{2})$.
	\item When a point crosses the $\theta_{0}=\pi/2$ plane we have that $\theta_{0}$ changes direction and $(\theta_{1},\theta_{2})$ is identified with $(\theta_{1},\theta_{2}+\pi)$.
\end{enumerate}

\mgap

The (round) Riemannian metric in Hopf coordinates is:
\begin{equation}\label{E:roundmetric}
	ds^{2}=d\theta_{0}^{2}+\sin^{2}\!\theta_{0}\,d\theta_{1}^{2}+\cos^{2}\!\theta_{0}\,d\theta_{2}^{2},
\end{equation}
and the Laplacian is
$$\Delta=\frac{1}{\sin(2\theta_{0})}\frac{\partial}{\partial\theta_{0}}\circ\sin(2\theta_{0})\frac{\partial}{\partial\theta_{0}}+\csc^{2}\!\theta_{0}\,\frac{\partial^{2}}{\partial\theta_{1}^{2}}+\sec^{2}\!\theta_{0}\,\frac{\partial^{2}}{\partial\theta_{2}^{2}}.$$

\mgap

We now write the sR structure in Hopf coordinates. We can introduce $r>0$, to give coordinates to $\Rbb^{4}$, allowing us to write the 
$\frac{\partial}{\partial x_{j}},\frac{\partial}{\partial y_{j}}$ in terms of the
$\frac{\partial}{\partial\theta_{j}},\frac{\partial}{\partial r}$. Then restricting to functions on $S^{3}$ we get:
\begin{align*}
	\frac{\partial}{\partial x_{1}}&=\cos\theta_{1}\cos\theta_{0}\frac{\partial}{\partial\theta_{0}}-\sin\theta_{1}\csc\theta_{0}\frac{\partial}{\partial\theta_{1}}\\
	\frac{\partial}{\partial y_{1}}&=\sin\theta_{1}\cos\theta_{0}\frac{\partial}{\partial\theta_{0}}+\cos\theta_{1}\csc\theta_{0}\frac{\partial}{\partial\theta_{1}}\\
	\frac{\partial}{\partial x_{2}}&=-\cos\theta_{2}\sin\theta_{0}\frac{\partial}{\partial\theta_{0}}-\sin\theta_{2}\sec\theta_{0}\frac{\partial}{\partial\theta_{2}}\\
	\frac{\partial}{\partial y_{2}}&=-\sin\theta_{2}\sin\theta_{0}\frac{\partial}{\partial\theta_{0}}+\cos\theta_{2}\sec\theta_{0}\frac{\partial}{\partial\theta_{2}}.
\end{align*}
Our vector fields are then
\begin{align*}
	V&=\frac{\partial}{\partial\theta_{1}}+\frac{\partial}{\partial\theta_{2}}\\
	E_{1}&=-\cos(\theta_{1}+\theta_{2})\frac{\partial}{\partial \theta_{0}}+\sin(\theta_{1}+\theta_{2})\cot\theta_{0}\frac{\partial}{\partial \theta_{1}}
		-\sin(\theta_{1}+\theta_{2})\tan\theta_{0}\frac{\partial}{\partial \theta_{2}}\\
	E_{2}&=-\sin(\theta_{1}+\theta_{2})\frac{\partial}{\partial \theta_{0}}-\cos(\theta_{1}+\theta_{2})\cot\theta_{0}\frac{\partial}{\partial \theta_{1}}
		+\cos(\theta_{1}+\theta_{2})\tan\theta_{0}\frac{\partial}{\partial \theta_{2}}.
\end{align*} 
The commutation relations hold, same as before, and the vector fields are still orthonormal (of course, with respect to the Riemannian metric in Hopf coordinates).

\mgap

The sR metric, written in Hopf coordinates, is
\begin{equation*}
	S=
	\begin{pmatrix}
		1&0&0\\
		0&\cos^{2}\!\theta_{0}\,\sin^{2}\!\theta_{0}&-\cos^{2}\!\theta_{0}\,\sin^{2}\!\theta_{0}\\
		0&-\cos^{2}\!\theta_{0}\,\sin^{2}\!\theta_{0}&\cos^{2}\!\theta_{0}\,\sin^{2}\!\theta_{0}
	\end{pmatrix}.
\end{equation*}
Indeed it is easy to check that $E_{1}$ and $E_{2}$ are orthonormal with respect to $S$, and $V$ is in the kernel of $S$. Written as a two-tensor, 
\begin{equation*}
	S=d\theta_{0}\otimes d\theta_{0}+\cos^{2}\!\theta_{0}\,\sin^{2}\!\theta_{0}\,(d\theta_{1}-d\theta_{2})\otimes(d\theta_{1}-d\theta_{2}).
\end{equation*}

\mgap

The sR Laplacian, written in Hopf coordinates, is
$$\Delta_{sR}=E_{1}^{2}+E_{2}^{2}=\frac{1}{\sin(2\theta_{0})}\frac{\partial}{\partial\theta_{0}}\circ\sin(2\theta_{0})\frac{\partial}{\partial\theta_{0}}+\left(\cot\theta_{0}\frac{\partial}{\partial\theta_{1}}-\tan\theta_{0}\frac{\partial}{\partial\theta_{2}}\right)^{2}.$$

\mgap

We can consider the sR metric as being the limit of certain penalty metrics, where the $V$ direction is penalized by a factor $\lambda>1$. After simple linear algebra (multiplying the $V$ direction by $\lambda$, multiplying the other directions by $1$, and \emph{then} applying the Riemannian metric), the $\lambda$-penalty metric is given by the matrix
\begin{equation*}
	P_{\lambda}=
	\begin{pmatrix}
		1&&0&&0\\
		\\
		0&&(\lambda^{2}-1)\sin^{4}\!\theta_{0}+\sin^{2}\!\theta_{0}&&(\lambda^{2}-1)\cos^{2}\!\theta_{0}\,\sin^{2}\!\theta_{0}\\
		\\
		0&&(\lambda^{2}-1)\cos^{2}\!\theta_{0}\,\sin^{2}\!\theta_{0}&&(\lambda^{2}-1)\cos^{4}\!\theta_{0}+\cos^{2}\!\theta_{0}
	\end{pmatrix}.
\end{equation*}
Indeed, one can check that in fact $V$, $E_{1}$, and $E_{2}$ are orthogonal with respect to this metric, that $E_{1}$ and $E_{2}$ have length $1$, and that $V$ has length $\lambda$. We can easily compute
$$\det P_{\lambda}=\lambda^{2}\cos^{2}\!\theta_{0}\,\sin^{2}\!\theta_{0}$$
and
\begin{equation*}
	P_{\lambda}^{-1}=
	\begin{pmatrix}
		1&&0&&0\\
		\\
		0&&\cot^{2}\!\theta_{0}+\lambda^{-2}&&\lambda^{-2}-1\\
		\\
		0&&\lambda^{-2}-1&&\tan^{2}\!\theta_{0}+\lambda^{-2}
	\end{pmatrix}.
\end{equation*}
From this we find that the $\lambda$-penalty Laplacian on $S^{3}$ is:
\begin{equation*}
	\Delta_{\lambda}=\frac{\partial^{2}}{\partial\theta_{0}^{2}}+2\cot(2\theta_{0})\frac{\partial}{\partial\theta_{0}}
+\left(\cot\theta_{0}\frac{\partial}{\partial\theta_{1}}-\tan\theta_{0}\frac{\partial}{\partial\theta_{2}}\right)^{2}+\lambda^{-2}\left(\frac{\partial}{\partial\theta_{1}}+\frac{\partial}{\partial\theta_{2}}\right)^{2}.
\end{equation*}
That is,
$$\Delta_{\lambda}=E_{1}^{2}+E_{2}^{2}+\lambda^{-2}V^{2},$$
as might have been expected.

\mgap

Montgomery discovered an example in which geodesics with respect to the $\lambda$-penalty metric converge (as $\lambda\to\infty$) to sR geodesics that do \emph{not} solve the sR geodesic equations, in contrast to the Riemannian setting; that is, Montgomery discovered so-called \emph{singular geodesics} \cite{R:Montg94a}. For the case of $S^{3}$ (and more generally, in the contact case), singular geodesics do not exist, so it suffices to study the geodesic equations, or, equivalently, Hamilton's equations \cite{R:Montg}.

\mgap

We denote the dual variable to $\theta_{j}$ by $\xi_{j}$. The sR Hamiltonian is then
\begin{equation*}
	H(\theta,\xi)
	=\tfrac{1}{2}\xi_{0}^{2}+\tfrac{1}{2}(\cot\!\theta_{0}\,\xi_{1}-\tan\!\theta_{0}\,\xi_{2})^{2}.
\end{equation*}
Hamilton's equations, giving the sR geodesics, are then, for $j=0,1,2$,
\begin{align*}
	\dot\theta_{j}&=\frac{\partial H}{\partial \xi_{j}}\\
	\dot\xi_{j}&=-\frac{\partial H}{\partial \theta_{j}}.
\end{align*}
Explicitly,
\begin{equation}\label{E:HamEqHopf}
\begin{aligned}
	\dot\theta_{0}&=\xi_{0}\\
	\dot\theta_{1}&=\cot^{2}\!\theta_{0}\,\xi_{1}-\xi_{2}\\
	\dot\theta_{2}&=\tan^{2}\!\theta_{0}\,\xi_{2}-\xi_{1}\\
	\dot\xi_{0}&=\cot\!\theta_{0}\,\csc^{2}\!\theta_{0}\,\xi_{1}^{2}-\tan\!\theta_{0}\,\sec^{2}\!\theta_{0}\,\xi_{2}^{2}\\
	\dot\xi_{1}&=0\\
	\dot\xi_{2}&=0.
\end{aligned}
\end{equation}
One obvious advantage of using Hopf coordinates is that $\xi_{1}$ and $\xi_{2}$ are constant along the flow; in addition, as always $H$ is constant along the flow, so we already have three conserved quantities. Also, these equations have a clear symmetry; for example, 
$$\cot(\tfrac{\pi}{2}-\theta_{0})\,\csc^{2}(\tfrac{\pi}{2}-\theta_{0})=\tan\theta_{0}\,\sec^{2}\theta_{0}.$$

\mgap

The penalty Hamiltonian is:
\begin{equation}\label{E:penHam}
	\begin{aligned}
	H_{\lambda}(\theta,\xi)
	&=H+\tfrac{1}{2\lambda^{2}}(\xi_{1}+\xi_{2})^{2}\\
	&=\tfrac{1}{2}\xi_{0}^{2}+\tfrac{1}{2}(\cot\!\theta_{0}\,\xi_{1}-\tan\!\theta_{0}\,\xi_{2})^{2}+\tfrac{1}{2\lambda^{2}}(\xi_{1}+\xi_{2})^{2}.
	\end{aligned}
\end{equation}
The corresponding penalty Hamiltonian equations, giving the penalty geodesics, are then
\begin{equation}\label{E:penHamEqHopf}
\begin{aligned}
	\dot\theta_{0}&=\xi_{0}\\
	\dot\theta_{1}&=\cot^{2}\!\theta_{0}\,\xi_{1}-\xi_{2}+\lambda^{-2}(\xi_{1}+\xi_{2})\\
	\dot\theta_{2}&=\tan^{2}\!\theta_{0}\,\xi_{2}-\xi_{1}+\lambda^{-2}(\xi_{1}+\xi_{2})\\
	\dot\xi_{0}&=\cot\!\theta_{0}\,\csc^{2}\!\theta_{0}\,\xi_{1}^{2}-\tan\!\theta_{0}\,\sec^{2}\!\theta_{0}\,\xi_{2}^{2}\\
	\dot\xi_{1}&=0\\
	\dot\xi_{2}&=0.
\end{aligned}
\end{equation}
For the case of the Riemannian metric on $S^{3}$, that is, the case $\lambda=1$, the equations simplify, and we get
\begin{equation*}
\begin{aligned}
	\dot\theta_{1}&=\csc^{2}\!\theta_{0}\,\xi_{1}\\
	\dot\theta_{2}&=\sec^{2}\!\theta_{0}\,\xi_{2}.
\end{aligned}
\end{equation*}
When $\lambda=1$, the solutions of Hamilton's equations are great circles on $S^{3}$.

\mgap

\section{Categorizing sR Geodesics}\label{S:Categ}

Our categorization of sR geodesics is based on a reduced problem. In Hamilton's equations (\ref{E:HamEqHopf}), since $\xi_{1}$ and $\xi_{2}$ are constant along the flow, we can isolate the equations
\begin{equation*}
\begin{aligned}
	\dot\theta_{0}&=\xi_{0}\\
	\dot\xi_{0}&=\cot\!\theta_{0}\,\csc^{2}\!\theta_{0}\,\xi_{1}^{2}-\tan\!\theta_{0}\,\sec^{2}\!\theta_{0}\,\xi_{2}^{2}\,\,\,,
\end{aligned}
\end{equation*}
which are Hamilton's equations for the sR Hamiltonian $H$ considered as a function of \emph{two} variables
\begin{equation}\label{E:onedimSchr}
	H(\theta_{0},\xi_{0})=\tfrac{1}{2}\xi_{0}^{2}+\tfrac{1}{2}(\cot\!\theta_{0}\,\xi_{1}-\tan\!\theta_{0}\,\xi_{2})^{2}.
\end{equation}
Equation (\ref{E:onedimSchr}) can be viewed as a \emph{one-dimensional} energy equation: it is of the form
$$\text{Energy}=\text{Kinetic Energy}+\text{Potential Energy},$$
with potential function
$$U=\tfrac{1}{2}(\cot\!\theta_{0}\,\xi_{1}-\tan\!\theta_{0}\,\xi_{2})^{2}.$$

\mgap

We now list the various disjoint cases:

\mgap

\begin{enumerate}
	\item[\textbf{1.}] A fixed point in the $(\theta_{0},\xi_{0})$ phase plane. From our original choice of coordinates we may assume that $\theta_{0}\equiv \pi/4$, and then $\xi_{1}^{2}=\xi_{2}^{2}$. 
	
	\begin{enumerate}
		\item[\textbf{1a.}] $\xi_{1}=\xi_{2}$. (This is precisely the case when $H=0$.) Then from Hamilton's equations $\theta_{0}$, 
			$\theta_{1}$, and $\theta_{2}$ are constant; this gives a degenerate sR geodesic of length $0$.
		\item[\textbf{1b.}] $\xi_{1}=-\xi_{2}\neq 0$. Hamilton's equations then say that the speed on the Hopf cube is $\sqrt{2}\,|\xi_{1}-\xi_{2}|$, and the length of the (simple) closed curve on the Hopf cube is $\sqrt{2}\,2\pi$, so the period is $2\pi/|\xi_{1}-\xi_{2}|$. On $S^{3}$ the speed is $|\xi_{1}-\xi_{2}|$, so the length of this closed sR geodesic is $2\pi$. See Figure~\ref{F:case1b}.
	\end{enumerate}
\end{enumerate}

\mgap

We categorize the remaining cases in terms of the potential function $U$.

\mgap

\begin{enumerate}
	\item[\textbf{2.}] The ``free'' case $U\equiv 0$. This happens precisely when $\xi_{1}=\xi_{2}=0$. (We have already dispensed with the case when $\theta_{0}$ is constant.) By Hamilton's equations $\dot\theta_{1}$, $\dot\theta_{2}$, and $\dot\xi_{0}$, are also identically zero, while $\dot\theta_{0}=\xi_{0}$. That is, we have a point with speed $|\xi_{0}|$ moving purely in the $\theta_{0}$-direction; the length of this (simple) closed geodesic is $2\pi$. (It is both a geodesic and a sR geodesic.) See Figure~\ref{F:case2}.
	
	\begin{figure}
	\centering
	\begin{minipage}{0.45\textwidth}
	\centering
	\includegraphics[width=0.8\textwidth]{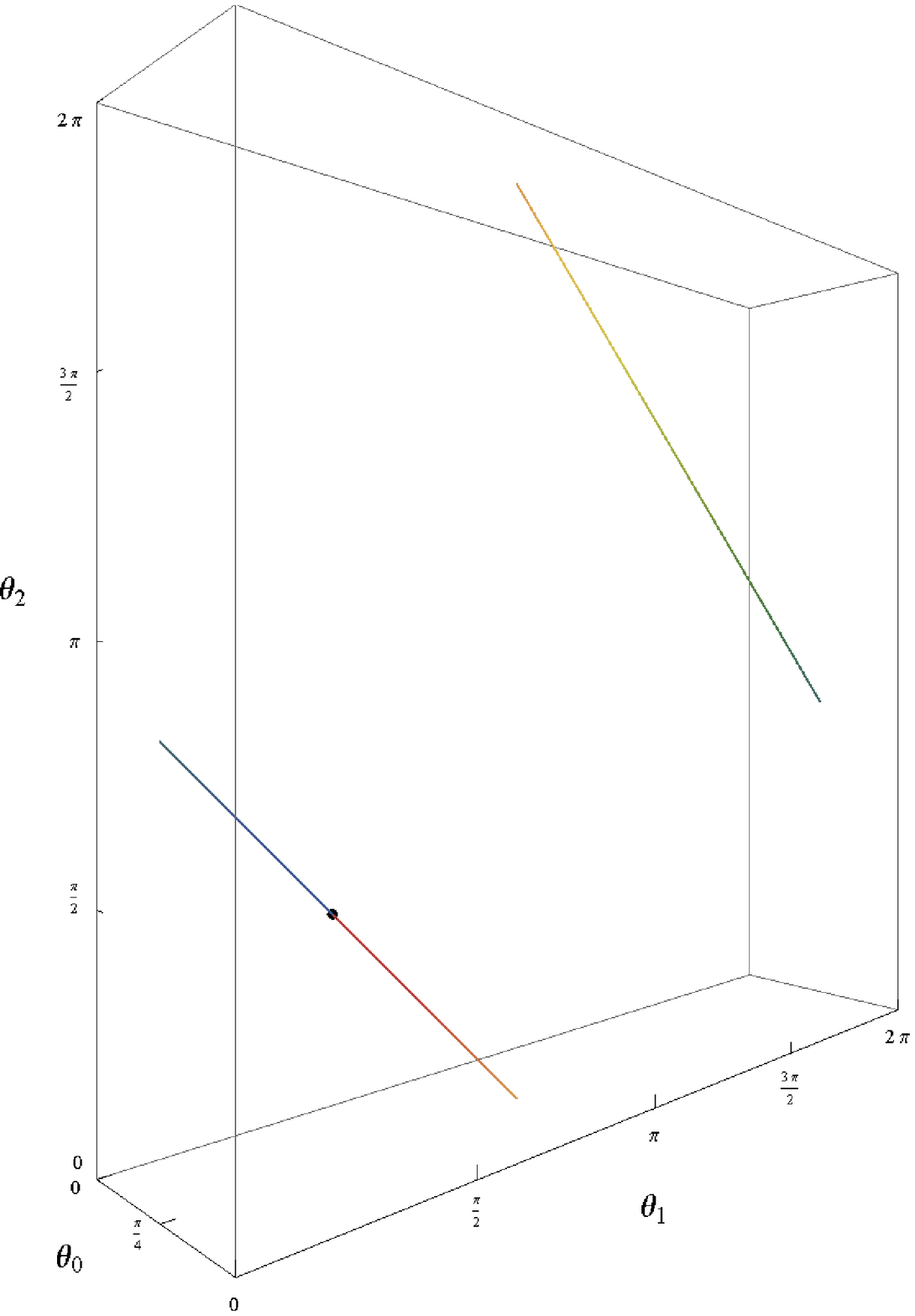} 
	\caption{Case 1b.}
	\label{F:case1b}
	\end{minipage}\hfill
	\begin{minipage}{0.45\textwidth}
	\centering
	\includegraphics[width=0.8\textwidth]{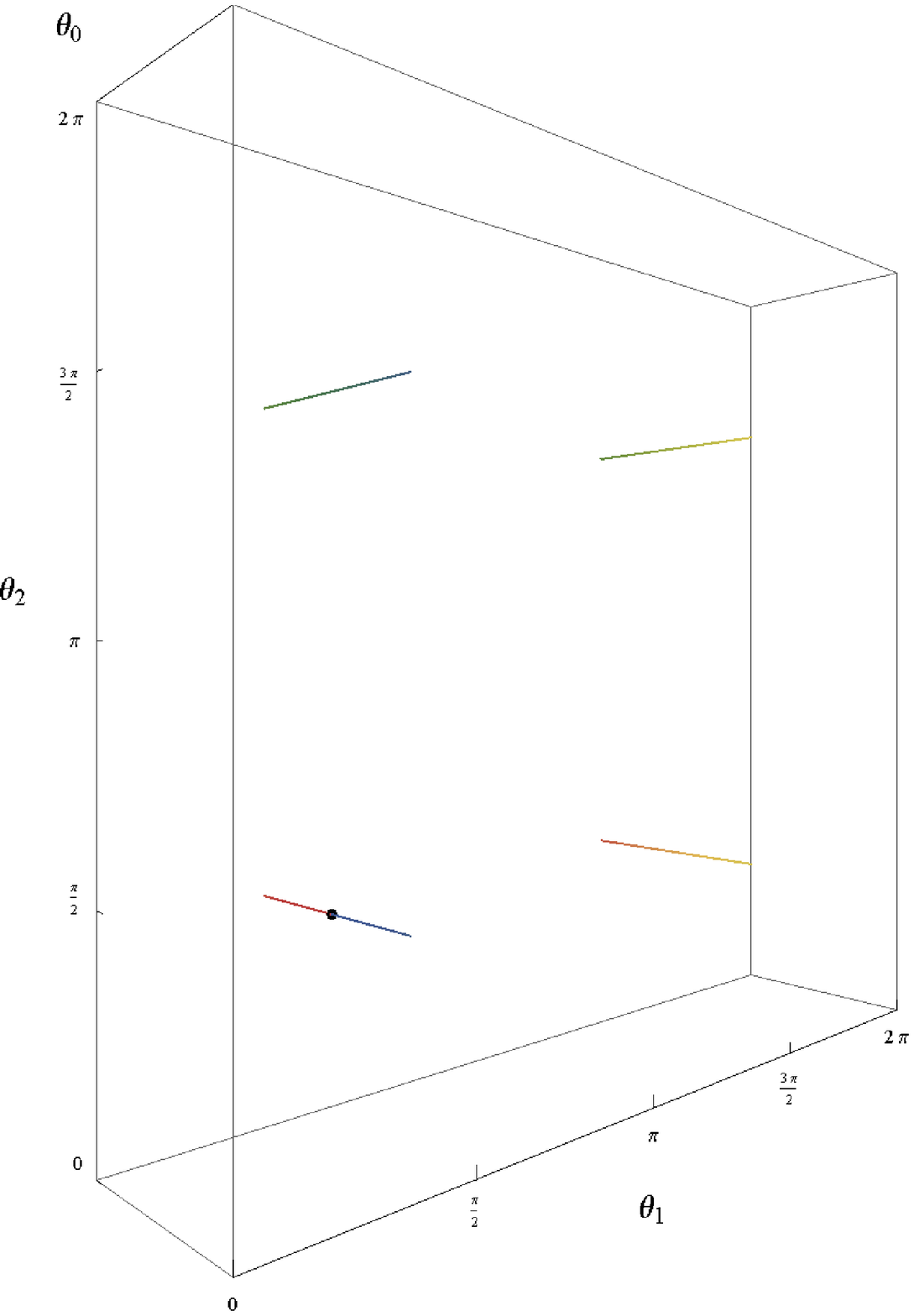}
	\caption{Case 2.}
	\label{F:case2}
	\end{minipage}
	\end{figure}

	\mgap

	\item[\textbf{3.}] $\xi_{1}\neq 0$ and $\xi_{2}\neq 0$. Then $U$ is a potential well with a single non-degenerate minimum occurring when $\tan^{4}\!\theta_{0}=\xi_{1}^{2}/\xi_{2}^{2}$. Typical potential functions are in Figures~\ref{F:potentialUss} and~\ref{F:potentialUos}, for $\xi_{1}$ and $\xi_{2}$ with the same and opposite signs, respectively.
	
	\sgap

	\begin{figure}
	\centering
	\begin{minipage}{0.45\textwidth}
	\centering
	\includegraphics[width=0.8\textwidth]{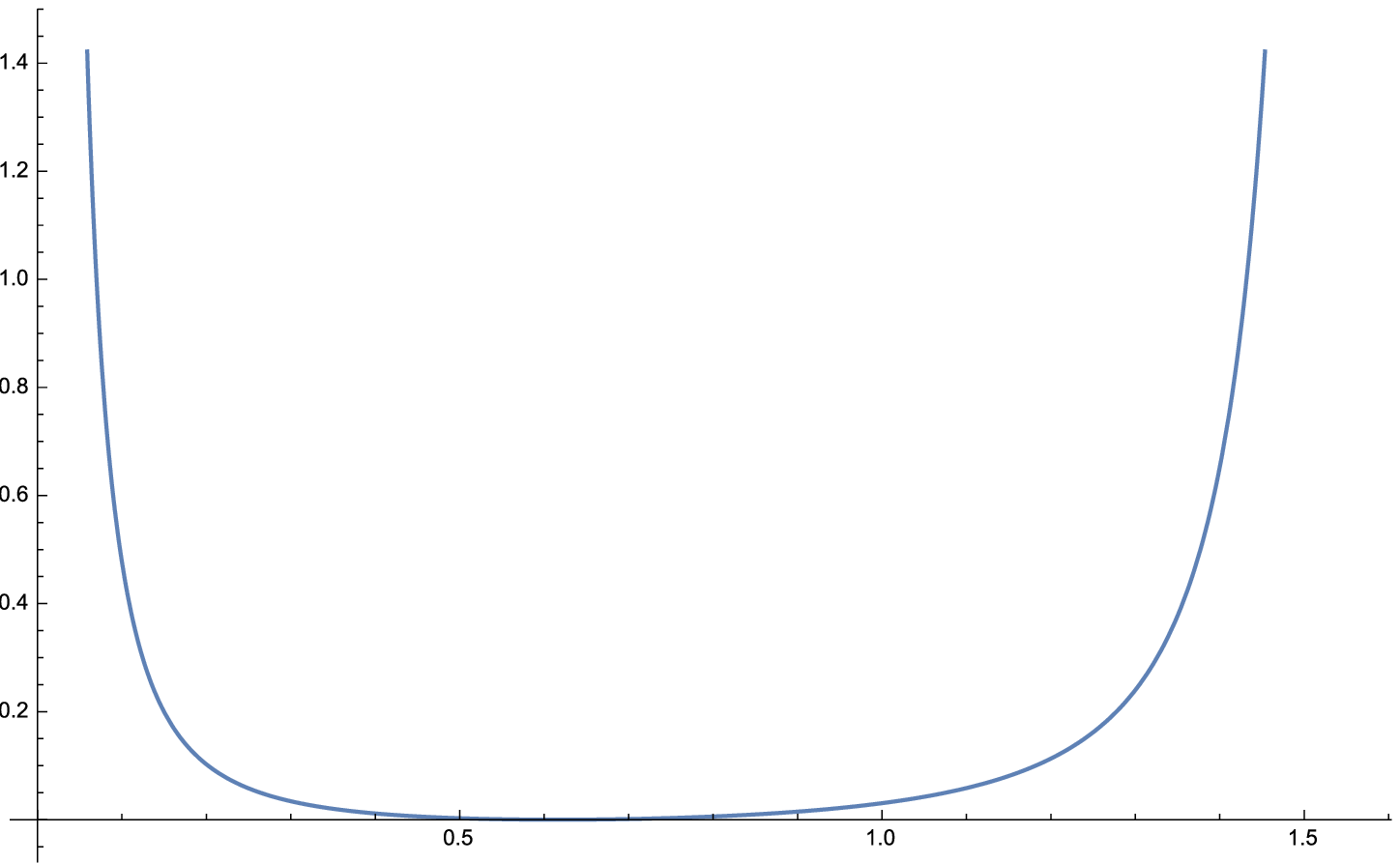}
	\caption{$\xi_{1}=0.1$, $\xi_{2}=0.2$.}
	\label{F:potentialUss}
	\end{minipage}\hfill
	\begin{minipage}{0.45\textwidth}
	\centering
	\includegraphics[width=0.8\textwidth]{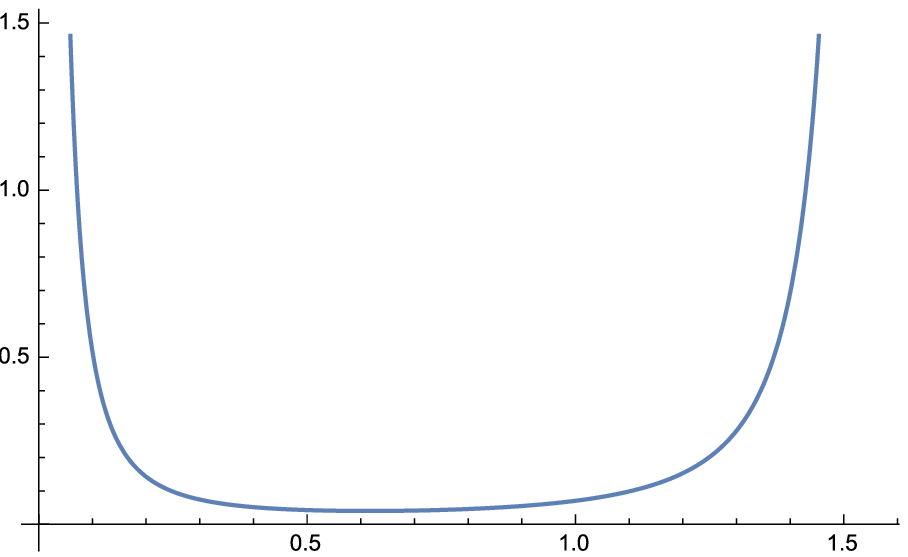} 
	\caption{$\xi_{1}=0.1$, $\xi_{2}=-0.2$.}
	\label{F:potentialUos}
	\end{minipage}
	\end{figure}
	
	\begin{figure}
	\centering
	\begin{minipage}{0.45\textwidth}
	\centering
	\includegraphics[width=0.8\textwidth]{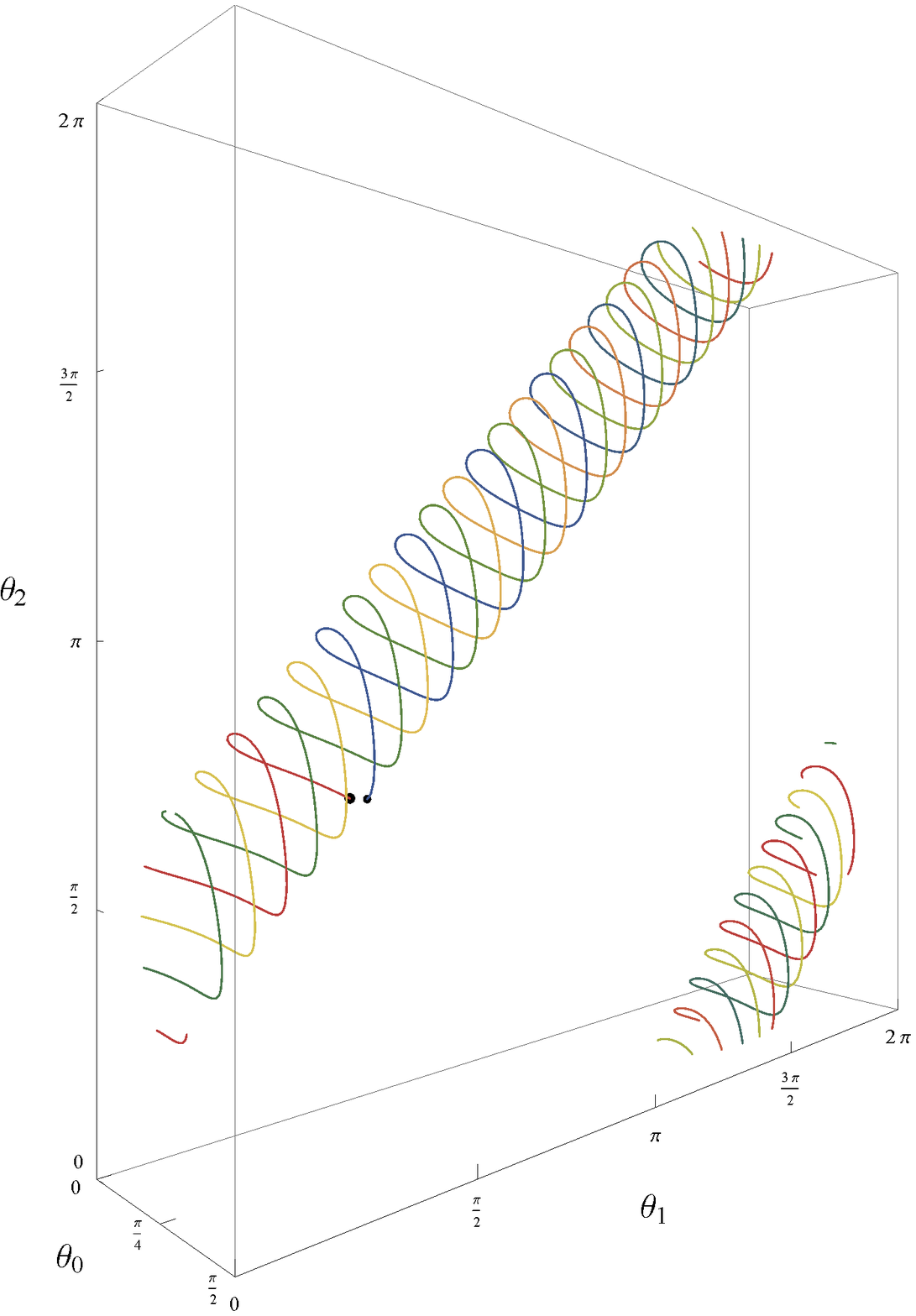} 
	\caption{Case 3. Ex. 1.}
	\label{F:helix1}
	\end{minipage}\hfill
	\begin{minipage}{0.45\textwidth}
	\centering
	\includegraphics[width=0.8\textwidth]{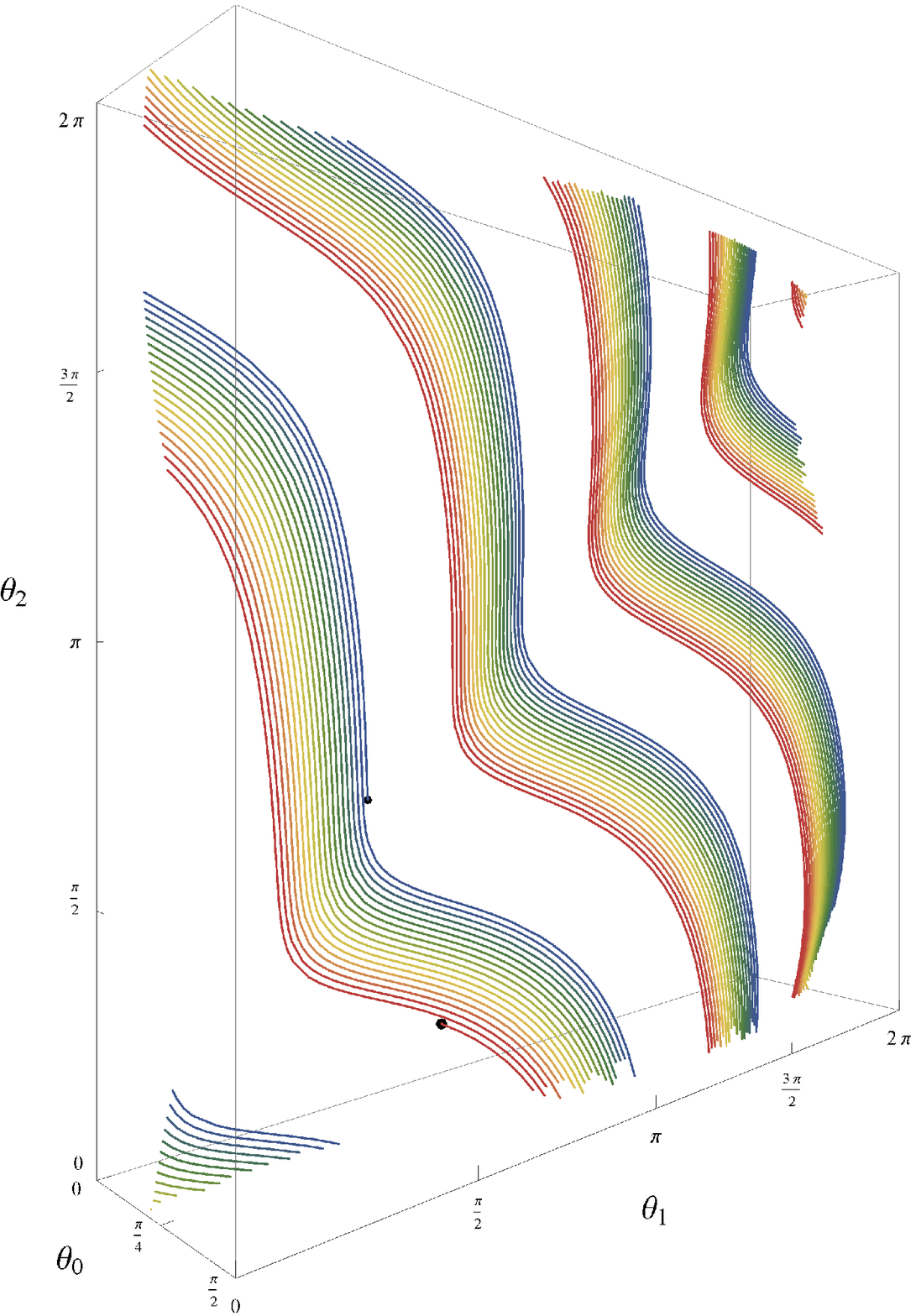} 
	\caption{Case 3. Ex. 2.}
	\label{F:blanket1}
	\end{minipage}
	\end{figure}

\sgap

Since in this case $\theta_{0}$ is \emph{not} constant, its period is
\begin{align*}
	\text{Period}(\theta_{0})
	&=2\int_{a}^{b}\frac{d\theta_{0}}{\sqrt{2(H-U)}}\\
	&=2\int_{a}^{b}\frac{d\theta_{0}}{\sqrt{2H-(\cot\!\theta_{0}\,\xi_{1}-\tan\!\theta_{0}\,\xi_{2})^{2}}}.
\end{align*}
Here $a$ and $b$ are the ``turning points,'' where the kinetic energy is zero. 

\mgap

Fortunately it is possible to evaluate this integral using freshman calculus. Substituting 
$$x=\cos^{2}\!\theta_{0}\,,\qquad 0<\theta_{0}<\pi/2,$$
we get
\begin{equation*}
	\text{Period}(\theta_{0})
	=
	\int_{\cos^{2}(b)}^{\cos^{2}(a)}\frac{dx}{\sqrt{[-2H-(\xi_{1}+\xi_{2})^{2}]x^{2}+2(H+\xi_{1}\xi_{2}+\xi_{2}^{2})x-\xi_{2}^{2}}}.
\end{equation*}
The limits of integration are exactly the points where the denominator vanishes (where the velocity is zero), and we recall that the Hamiltonian for Riemannian geodesics is $H_{1}=H+\tfrac{1}{2}(\xi_{1}+\xi_{2})^{2}$ (the case $\lambda=1$), so we have
\begin{equation*}
	\text{Period}(\theta_{0})=\frac{1}{\sqrt{2H_{1}}}\int_{\cos^{2}(b)}^{\cos^{2}(a)}
	\frac{dx}{\sqrt{(\cos^{2}(a)-x)(x-\cos^{2}(b))}}.
\end{equation*}
This is an integral known to be solvable by elementary functions. Following page 366 of \emph{Advanced Calculus (New Edition)} by Frederick S. Woods\footnote{This is the book Feynman mentions in \emph{Surely You're Joking, Mr. Feynman!} as giving him valuable tricks for integration \cite{R:MrFeynman}.} \cite{R:Woods}, we make the substitution defined by
$$z^{2}+1=\frac{\cos^{2}(a)-\cos^{2}(b)}{x-\cos^{2}(b)}$$ and finally get the answer
\begin{align*}
	\text{Period}(\theta_{0})
	&=\sqrt{\frac{2}{H_{1}}}\int_{0}^{\infty}\frac{dz}{z^{2}+1}\\
	&=\frac{\pi}{\sqrt{2H_{1}}}.
\end{align*}
For future reference, we note that this is one-half the period of the (Riemannian) geodesic flow; after all, the speed of the geodesic flow is $\sqrt{2H_{1}}$, and we know the length of each geodesic, a great circle on $S^{3}$, to be $2\pi$. (See Section~\ref{S:closure}.)

\mgap

Examples are pictured in Figures~\ref{F:helix1} and~\ref{F:blanket1}, where $\xi_{1}$ and $\xi_{2}$ have the same and opposite signs, respectively.

\mgap

	\item[\textbf{4}.] It remains to check the exceptional cases when \{$\xi_{1}=0$ and $\xi_{2}\neq 0$\} and when \{$\xi_{1}\neq 0$ and $\xi_{2}=0$\}. For example, when $\xi_{1}=0$ the potential function is
$$U=\tfrac{1}{2}\tan^{2}\!\theta_{0}\,\xi_{2}^{2},\qquad 0<\theta_{0}<\pi/2.$$
The force induced by this potential causes the point to exit the Hopf cube through the $\theta_{0}=0$ plane; rather we interpret it as bouncing off the plane, returning to the Hopf cube but with $\theta_{1}$ shifted by $\pi$. (See Section~\ref{S:Hopf}.) With reasoning as in the previous case, we find that again $\text{Period}(\theta_{0})=\pi/\sqrt{2H_{1}}$. The case when $\xi_{1}\neq 0$ and $\xi_{2}=0$ follows by renaming the variables $\theta_{0}\leftrightarrow \frac{\pi}{2}-\theta_{0}$ and $\xi_{1}\leftrightarrow\xi_{2}$. In Figures~\ref{F:case4a} and~\ref{F:case4b} we have the cases when $0<\xi_{1}\ll \xi_{2}$ and $0<\xi_{2}\ll \xi_{1}$, respectively, which illustrate how exiting the Hopf cube and re-entering after a $\pi$-shift appears as a limiting case.

	\begin{figure}
	\centering
	\begin{minipage}{0.45\textwidth}
	\centering
	\includegraphics[width=0.8\textwidth]{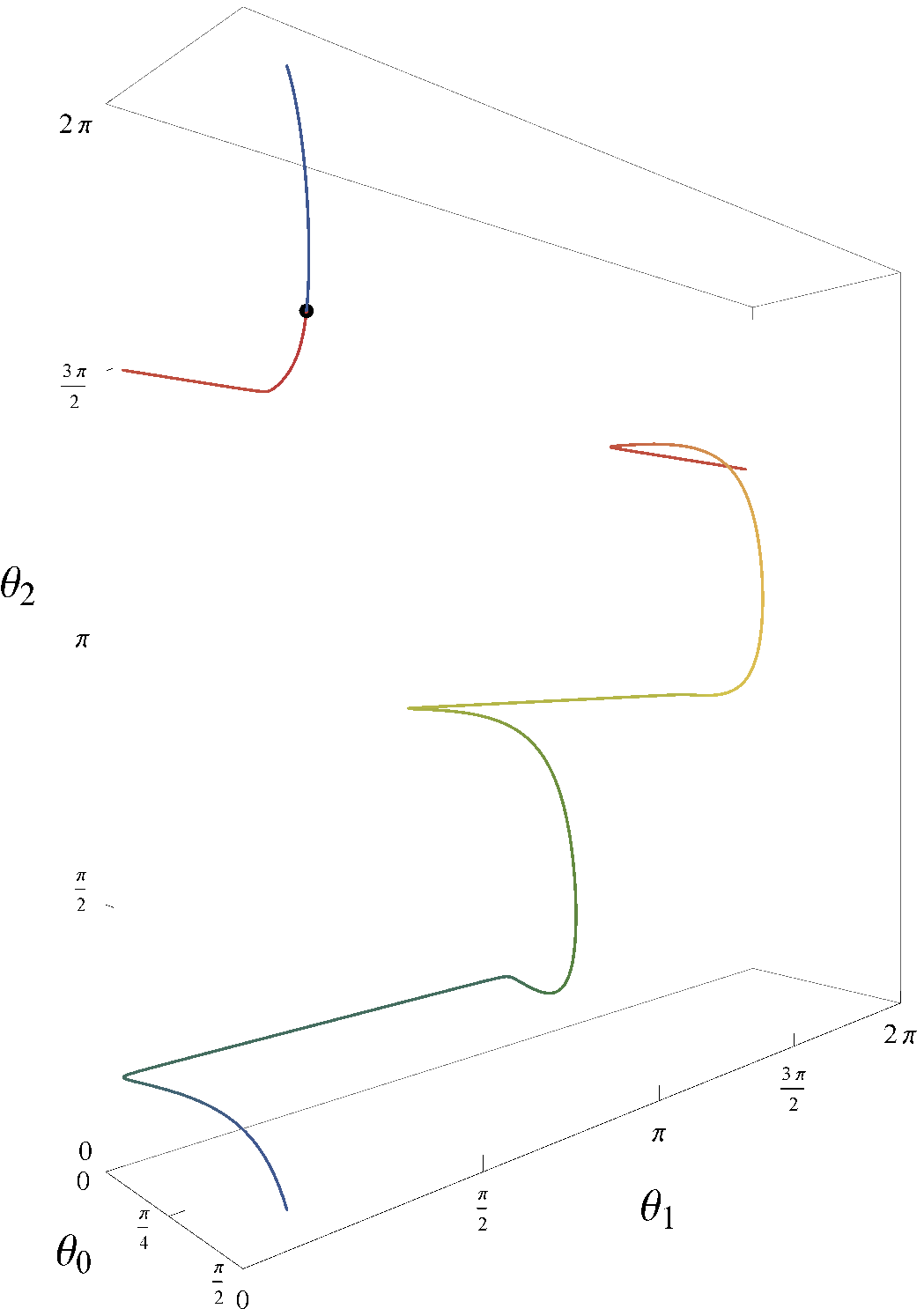} 
	\caption{$0<\xi_{1}\ll \xi_{2}$.}
	\label{F:case4a}
	\end{minipage}\hfill
	\begin{minipage}{0.45\textwidth}
	\centering
	\includegraphics[width=0.8\textwidth]{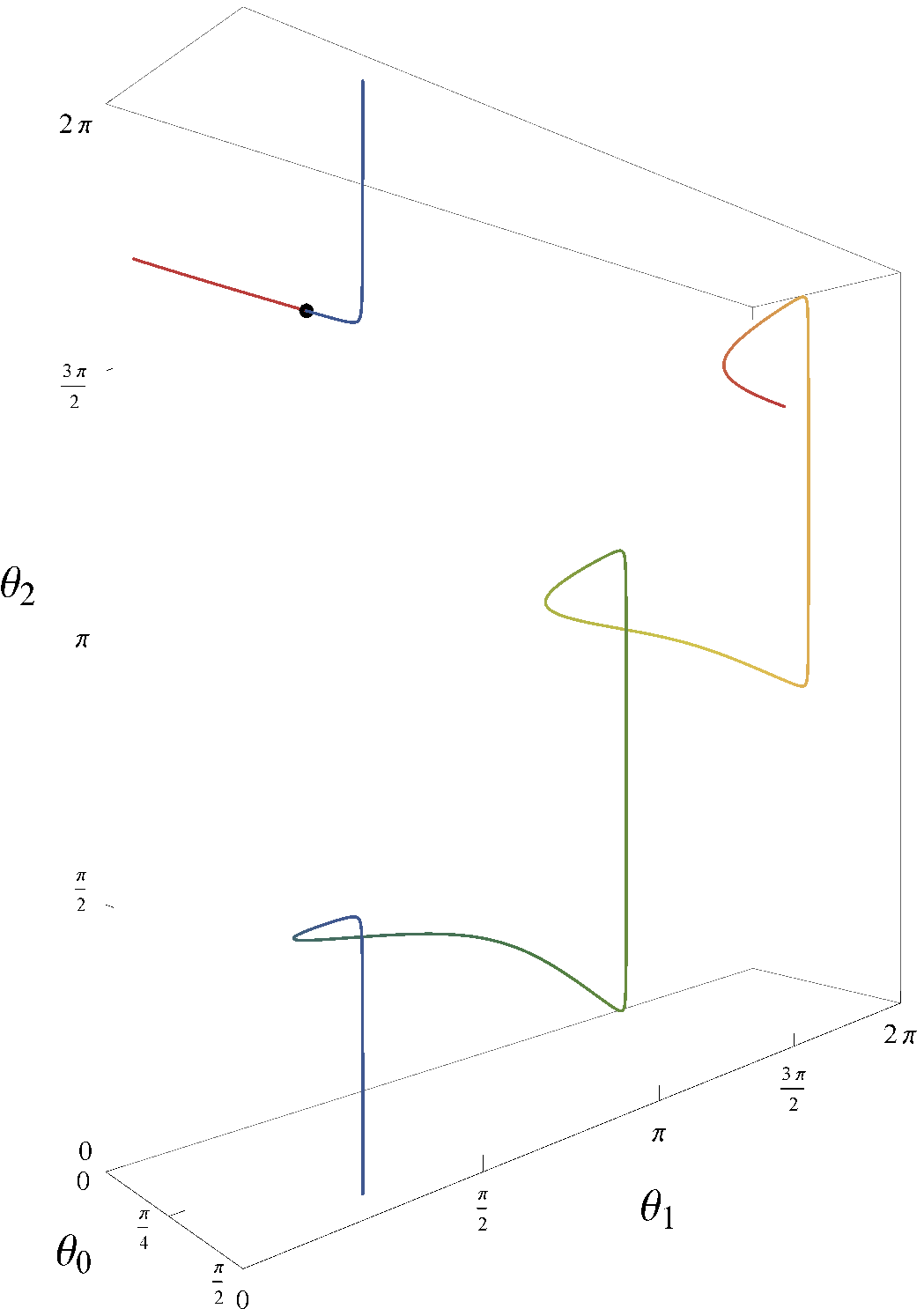} 
	\caption{$0<\xi_{2}\ll \xi_{1}$.}
	\label{F:case4b}
	\end{minipage}
	\end{figure}

\end{enumerate}

\sgap

Finally, we note that all sR geodesics are simple curves; that is, they do not self-intersect except trivially for closed curves. In Cases~1 and~2 above it is obvious. In Cases~3 and~4 we only need to wait until $\xi_{0}$ is zero, corresponding to the $\theta_{0}$-particle having zero kinetic energy in the potential well $U$. When $(\theta_{0},\theta_{1},\theta_{2})$ returns to that value, clearly $\xi_{0}$ is zero again, $\xi_{1},\xi_{2}$ are the same as always, and the $\dot\theta_{j}$ and $\dot\xi_{j}$ return to their values; thus the curve only self-intersects in the case of a closed curve, at the end of a period.

\sgap

\section{Determining Which sR Geodesics are Closed}\label{S:closure}

In this section we identify the closed sR geodesics on $S^{3}$; we only need to consider the Cases~3 and~4, and we may assume that the initial value of $\xi_{0}$ is zero. (See the comment at the end of Section~\ref{S:Categ}.) Hurtado and Rosales \cite{R:HurtadoRosales} found a necessary and sufficient condition in terms of geodesic curvature (see also the expository paper of D'Angelo and Tyson \cite{R:DAngTyson}):

\begin{theorem*} \cite{R:HurtadoRosales}
	Let $\gamma:\,\,\Rbb\to S^{3}$ be a complete sR geodesic of curvature $\lambda$. Then $\gamma$ is a closed curve diffeomorphic to a circle if and only if $\lambda/\sqrt{1+\lambda^{2}}$ is a rational number. Otherwise $\gamma$ is diffeomorphic to $\Rbb$ and is dense in some group translate of a Clifford torus.
\end{theorem*}

\mgap

Their proof relies on closed-form expressions of the sR geodesics. Here we give a condition which does not rely on closed-form expressions.

\mgap

From the $\lambda$-penalty Hamilton's equations (\ref{E:penHamEqHopf}), we see that the sR Hamiltonian vector field for the Hamiltonian $H$ is the difference of the Hamiltonian vector fields for the Hamiltonians $H_{1}$ and $H_{V}=\tfrac{1}{2}(\xi_{1}+\xi_{2})^{2}$. Moreover, the vector fields Lie-commute (it is easy to see that the Poisson bracket of $H_{1}$ and $H_{V}$ is zero), so the Hamiltonian flows for $H_{1}$ and $H_{V}$ commute. We can thus consider the $H$-flow as an $H_{1}$-flow followed by an $H_{V}$-flow.

\mgap

The Hamiltonian for the Riemannian geodesics may be written as
\begin{equation*}
	H_{1}(\theta,\xi)
	=\tfrac{1}{2}\xi_{0}^{2}+\tfrac{1}{2}(\csc^{2}\!\theta_{0}\,\xi_{1}^{2}+\sec^{2}\!\theta_{0}\,\xi_{2}^{2})\,\,,
\end{equation*}
(the penalty Hamiltonian (\ref{E:penHam}) with $\lambda=1$), so the first of Hamilton's equations, giving the velocities, are then
\begin{equation*}	
	\begin{aligned}	
		\dot\theta_{0}&=\xi_{0}\\
		\dot\theta_{1}&=\csc^{2}\!\theta_{0}\,\xi_{1}\\
		\dot\theta_{2}&=\sec^{2}\!\theta_{0}\,\xi_{2}.
	\end{aligned}
\end{equation*}
We see that the speed (measured using the Riemannian metric (\ref{E:roundmetric})) is $\sqrt{2H_{1}}$, which is constant. Moreover, the length of the Riemannian geodesic is $2\pi$, being a great circle, so that the period of the closed orbit is $2\pi/\sqrt{2H_{1}}=2\times\text{Period}(\theta_{0})$.

\mgap

On the other hand, the Hamiltonian $H_{V}$ has Hamiltonian equations
\begin{equation*}	
	\begin{aligned}	
		\dot\theta_{1}&=\xi_{1}+\xi_{2}\\
		\dot\theta_{2}&=\xi_{1}+\xi_{2}\\
		\dot\xi_{1}&=\dot\xi_{2}=0.
	\end{aligned}
\end{equation*}
Thus the speed (with respect to the Euclidean metric on the Hopf cube) is $\sqrt{2}|\xi_{1}+\xi_{2}|$. The length of the orbit (a circle fiber of the Hopf fibration) is $\sqrt{2}\cdot 2\pi$, so the period of the $H_{V}$-flow is $2\pi/|\xi_{1}+\xi_{2}|$. [It might seem strange that we find the speed and length with respect to the \emph{Euclidean} metric on the Hopf cube, but the Euclidean metric is sufficient to compute the period of the $H_{V}$-flow.]

\mgap

For a combination of an $H_{1}$-flow and an $H_{V}$-flow to result in a closed curve, we need the $H_{1}$-flow to return $\theta_{0}$ to its original value (since the $H_{V}$-flow has no $\frac{\partial}{\partial\theta_{0}}$ component). Thus the time elapsed must be an integer multiple of $\text{Period}(\theta_{0})=\pi/\sqrt{2H_{1}}$. If the integer is odd, the $H_{1}$-flow takes the point to its antipodal point, and we would need a half-period of the $H_{V}$-flow to return to the starting point. If the integer is even, the $H_{1}$-flow takes the point back to itself, and we could only allow full periods of the $H_{V}$-flow. To summarize, a necessary and sufficient condition for a closed sR geodesic is:
\begin{equation*}
	\text{time elapsed}=p\times \frac{\pi}{|\xi_{1}+\xi_{2}|}=q\times\frac{\pi}{\sqrt{2H_{1}}},
\end{equation*}
where $p,q\in\{1,2,3,\ldots\}$ are either both odd or both even. In particular,
\begin{equation}\label{E:closedCondHH}
	\frac{p}{q}=\frac{|\xi_{1}+\xi_{2}|}{\sqrt{2H_{1}}}=\sqrt{1-\frac{H}{H_{1}}}\in\mathbb{Q}\cap (0,1),
\end{equation}
The quantity $p/q$ is conserved along the flow and is positively homogeneous of degree zero in the $\xi$-variables. The condition (\ref{E:closedCondHH}) is also sufficient to have a closed sR geodesic. If it holds, then we have
$$H\text{-Period}=p\times \frac{\pi}{|\xi_{1}+\xi_{2}|}=q\times\frac{\pi}{\sqrt{2H_{1}}},$$
for the least such integers $0<p< q$ that are either both odd or both even.

\mgap

\begin{figure}
\includegraphics[width=0.35\textwidth]{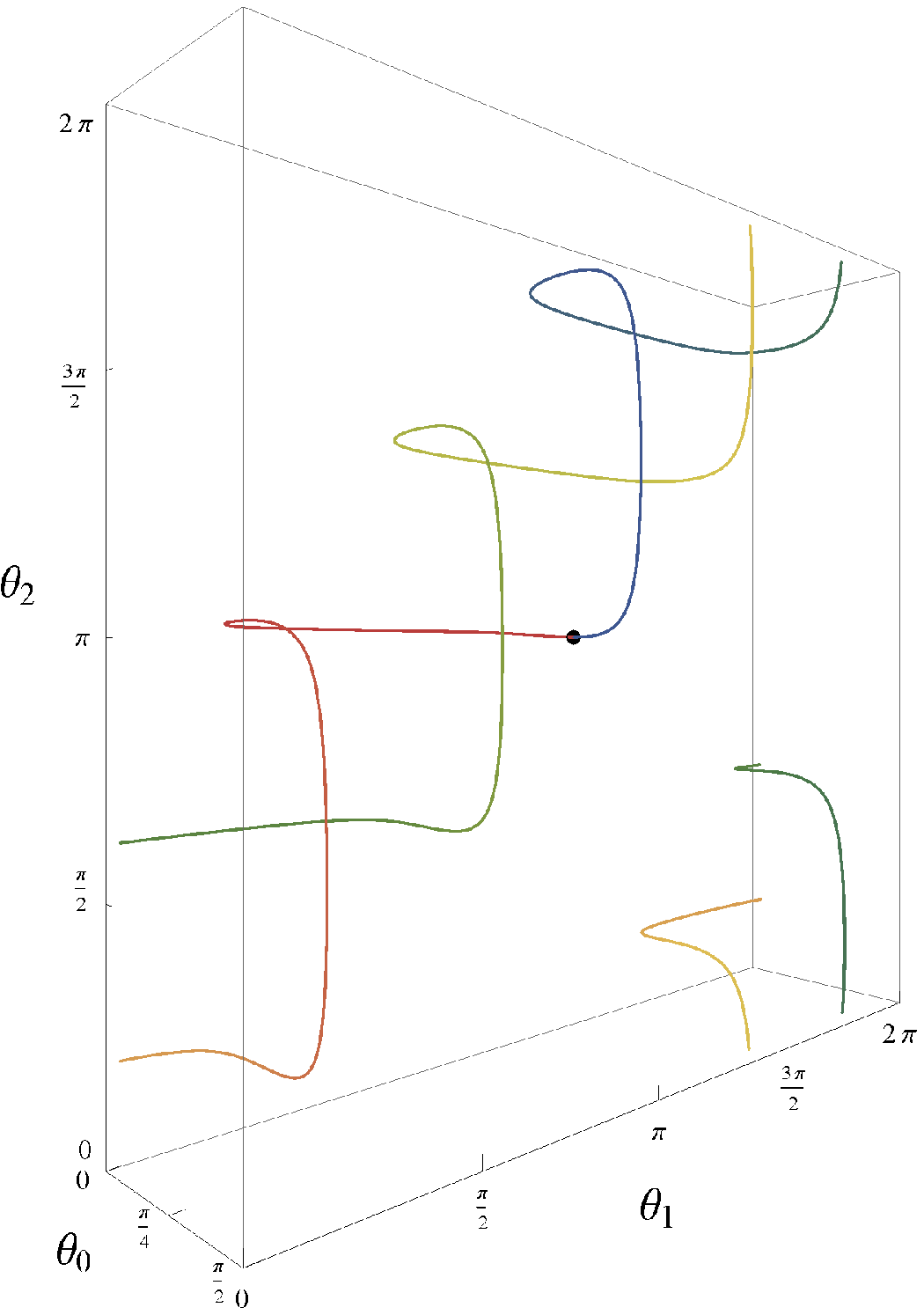}
\caption{An example with $r=1/5$.}\label{F:fiveLoops}
\end{figure}

When plotting sR geodesics in Cases~3 and~4, we can fix any $r\in\mathbb{Q}\cap(0,1)$ and re-write the closure condition (\ref{E:closedCondHH}) as
$$\xi_{0}^{2}=\frac{(\xi_{1}+\xi_{2})^{2}}{r^{2}}-\csc^{2}\!\theta_{0}\,\xi_{1}^{2}-\sec^{2}\!\theta_{0}\,\xi_{2}^{2}.$$
We can always find initial conditions satisfying this. Indeed, in Case~3 we can take any nonzero $\xi_{1}$ and $\xi_{2}$ and then take $\theta_{0}$ to maximize the right-hand side: $\tan^{2}\!\theta_{0}=|\xi_{1}/\xi_{2}|$. If $\xi_{1}$ and $\xi_{2}$ have the same sign, the right-hand side is always positive. If $\xi_{1}$ and $\xi_{2}$ have opposite signs, we need
$$\left|\frac{\xi_{1}+\xi_{2}}{\xi_{1}-\xi_{2}}\right|>r,$$
which is only valid for certain $\xi_{1}$ and $\xi_{2}$. Case~4 is similar. Then we can solve for $\xi_{0}$, use those numbers as the initial conditions in Hamilton's equations, and then plot the closed sR geodesic. Taking, for example, $r=1/5$, $\xi_{1}=0.6$, and $\xi_{2}=0.7$ we get the sR geodesic in Figure~\ref{F:fiveLoops}.

\mgap

\section{The sR Length Spectrum}\label{S:lengths}

To calculate the lengths of the closed sR geodesics we again only need to consider Cases~3 and~4 (the cases where $\theta_{0}$ oscillates). We found in the previous section that a sR geodesic is closed when the period of the $H_{1}$-flow and the period of the $H_{V}$-flow are commensurable. Then we have
\begin{equation}\label{E:periodH}
	\text{Period of $H$-flow}=p\times \frac{\pi}{|\xi_{1}+\xi_{2}|}=q\times\frac{\pi}{\sqrt{2H_{1}}},
\end{equation}
for the least such integers $0<p< q$ where $p,q$ are either both odd or both even. Since we know the speed of the sR geodesic is a constant $\sqrt{2H}$, we have that the length is
\begin{equation}\label{E:lengthPart1}
\begin{aligned}
	\text{Length}&=\text{Period}\times\text{Speed}\\
	&=\frac{\pi q}{\sqrt{2H_{1}}}\times\sqrt{2H}\\
	&=\pi q\sqrt{\frac{H}{H_{1}}}\\
	&=\pi\sqrt{q^{2}-p^{2}},
\end{aligned}
\end{equation}
for the least integers $0<p<q$ satisfying (\ref{E:periodH}) where $p,q$ are either both odd or both even.

\mgap

\begin{figure}
\includegraphics[width=0.4\textwidth]{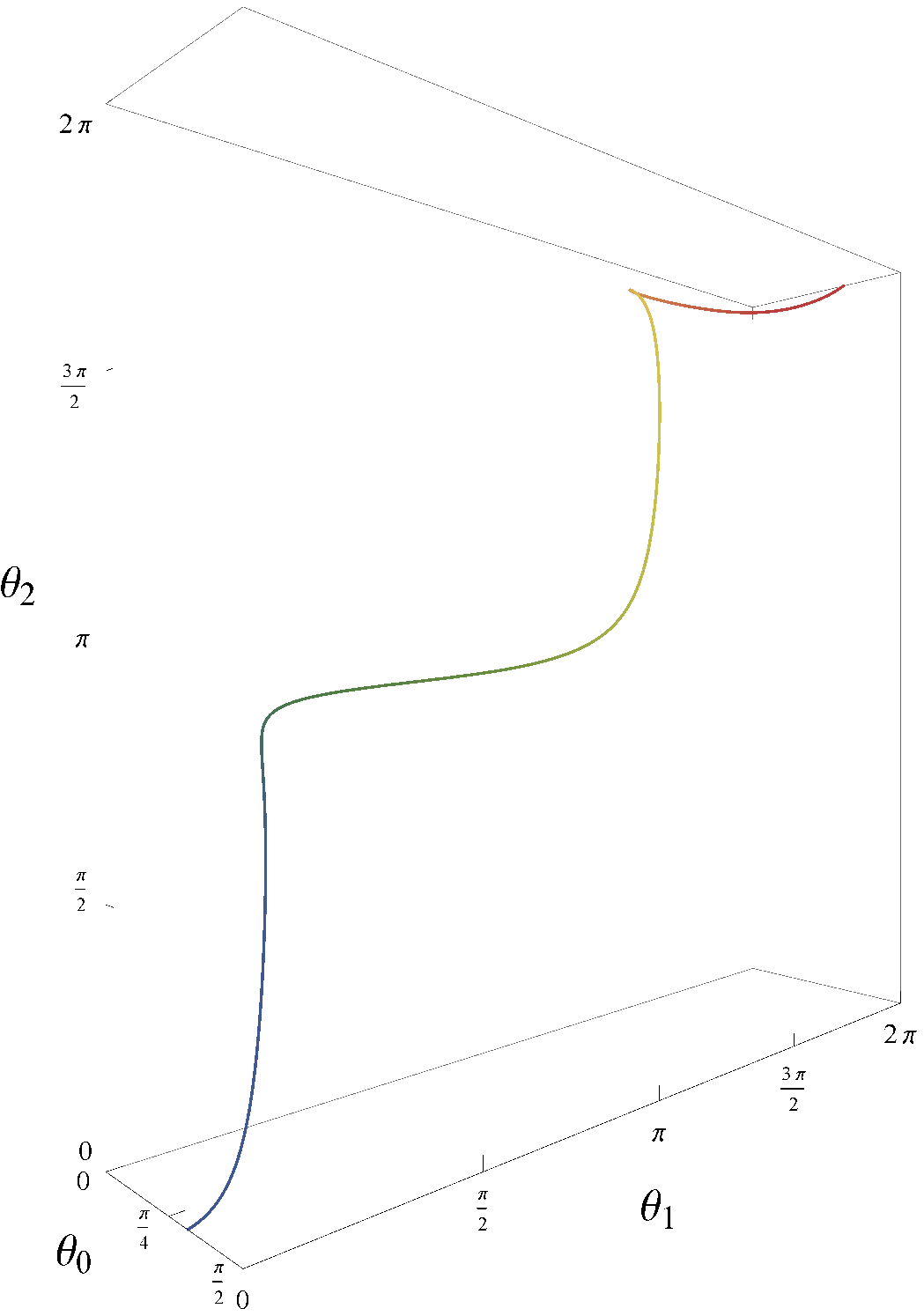}
\caption{A Riemannian geodesic in Hopf coordinates.}\label{F:RiemTwoLoops}
\end{figure}

We have another formulation of length that explains the repeating patterns seen in the figures. We know that the distance traveled in one $\theta_{0}$-period is 
\begin{equation*}
	\text{Period}(\theta_{0})\times \text{speed}=\pi\sqrt{\frac{H}{H_{1}}}.
\end{equation*}
Thus the length of a closed sR geodesic is
\begin{equation*}
	\text{Length}=\pi\times\text{(number of $\theta_{0}$-oscillations)}\times\sqrt{\frac{H}{H_{1}}}.
\end{equation*}
Comparing with equation (\ref{E:lengthPart1}), we find that
\begin{equation*}
	\text{number of $\theta_{0}$-oscillations}=q.
\end{equation*}
Moreover, we see from Hamilton's equations that the curve segments traced out by $\theta_{0}$-oscillations are congruent to each other. A similar argument shows that Riemannian geodesics in Hopf coordinates consist of two $\theta_{0}$-oscillations, as illustrated in Figure~\ref{F:RiemTwoLoops}.

\mgap

To summarize, we have found that if a sR geodesic is closed then the initial conditions must satisfy
$$\sqrt{1-\frac{H}{H_{1}}}=\frac{|\xi_{1}+\xi_{2}|}{\sqrt{2H_{1}}}\in\mathbb{Q}\cap(0,1)$$
and that the length of the closed sR geodesic is 
\begin{equation*}
	\text{Length}=\pi\sqrt{q^{2}-p^{2}},
\end{equation*}
for the least integers $0<p<q$ satisfying (\ref{E:periodH}) where $p,q$ are either both odd or both even.

\mgap

In fact, every such number is attained as a length; we simply follow the procedure:
\begin{enumerate}[(i)]
	\item Choose any $p/q\in\mathbb{Q}\cap(0,1)$, with $\gcd(p,q)=1$.
	\item As seen at the end of Section~\ref{S:closure}, we can choose initial conditions so that
	$$\frac{p}{q}=\sqrt{1-\frac{H}{H_{1}}}=\frac{|\xi_{1}+\xi_{2}|}{\sqrt{2H_{1}}}.$$
	Thus 
	$$p\times\frac{\pi}{|\xi_{1}+\xi_{2}|}=q\times\frac{\pi}{\sqrt{2H_{1}}}.$$
	\item If $p$ and $q$ are both odd, then the sR geodesic with those initial conditions has length $\pi\sqrt{q^{2}-p^{2}}.$
	If one of $\{p,q\}$ is odd and the other is even, the sR geodesic with those initial conditions has length $2\pi\sqrt{q^{2}-p^{2}}.$
\end{enumerate}

\mgap

 Thus the length spectrum consists of $2\pi$ and the numbers
 $$\pi\sqrt{q^{2}-p^{2}}$$
 where $0<p<q$ are odd integers with $\gcd(p,q)=1$, and 
$$2\pi\sqrt{q^{2}-p^{2}}$$
 where $0<p<q$ are integers, one odd and the other even, with $\gcd(p,q)=1$.

\mgap

We now give an alternative characterization of these numbers. It is simpler to work with squares of lengths divided by $\pi^{2}$. Then we wish to characterize the set $S$ of numbers consisting of $4$ and
$$\epsilon(q^{2}-p^{2}),$$
where $0<p<q$ are integers with $\gcd(p,q)=1$ and
\begin{equation*}
	\epsilon=
	\begin{cases}
		1&\text{if $p$ and $q$ are both odd}\\
		4&\text{if one of $\{p,q\}$ is odd and the other is even}.
	\end{cases}
\end{equation*}

\mgap

In the $\epsilon=1$ case we take the examples $p=2k-1$ and $q=2k+1$, $k\in\mathbb{N}$, to get
$$q^{2}-p^{2}=4(2k),\qquad k\in\mathbb{N}.$$
In the $\epsilon=4$ case, we take the examples $p=k$ and $q=k+1$, $k\in\mathbb{N}$, to get
$$4(q^{2}-p^{2})=4(2k+1),\qquad k\in\mathbb{N}.$$ This shows that $4\mathbb{N}\subset S$.
Now suppose that $n\in S$ and $4\nmid n$. Then clearly $n$ can only be in the $\epsilon=1$ case, so there would be odd integers $0<p<q$ with $\gcd(p,q)=1$ such that $n=q^{2}-p^{2}$. This is easily seen to be impossible. Thus in fact $4\mathbb{N}=S$.

\mgap

We note that if $n\in S$ and $8\mid n$, then $n$ cannot be in the $\epsilon=4$ case, and that if $n\in S$ and $n=4(2k+1)$, $k\in\mathbb{N}$, then $n$ cannot be in the $\epsilon=1$ case. Both of these statements easily follow from parity arguments.

\mgap

Converting back to the language of lengths, we find that the set of lengths of the closed sR geodesics is
$$\{2\pi \sqrt{n};\,\,n\in\mathbb{N}\}.$$
By the previous paragraph, odd $n$ correspond to ``full periods'' of the $H_{V}$-flow and geodesic flow (the $\epsilon=4$ case), and even $n$ correspond to ``half periods'' of both the $H_{V}$-flow and geodesic flow (the $\epsilon=1$ case).

\mgap

\section{The Spectrum of the subLaplacian}\label{S:spectrum}

The subLaplacian $-\Delta_{sR}$ has a compact resolvent, hence has a pure discrete spectrum $0<\lambda_{1}\leq \lambda_{2}\leq \cdots \leq\lambda_{n}\leq \cdots$, with $\lambda_{n}\to +\infty$ as $n\to +\infty$, and a complete orthonormal set of eigenfunctions. (See, for example, the recent paper of Colin de Verdi\`{e}re, et al. \cite{R:CdVproc}.) In fact, in the case of $S^{3}$, the eigenfunctions of the subLaplacian are the same as the eigenfunctions of the Laplacian. We recall that $\Delta_{\lambda}=E_{1}^{2}+E_{2}^{2}+\lambda^{-2}V^{2}$ is the $\lambda$-penalty Laplacian, with $\lambda=1$ giving the Riemannian Laplacian on the sphere $\Delta_{S^{3}}$, and $\lambda=\infty$ giving the subLaplacian on the sphere. In Hopf coordinates we have
\begin{equation*}
	\Delta_{sR}=E_{1}^{2}+E_{2}^{2}
	=\frac{1}{\sin(2\theta_{0})}\frac{\partial}{\partial\theta_{0}}\circ\sin(2\theta_{0})\frac{\partial}{\partial\theta_{0}}
		+\left(\cot\theta_{0}\frac{\partial}{\partial\theta_{1}}-\tan\theta_{0}\frac{\partial}{\partial\theta_{2}}\right)^{2}.
\end{equation*}
It is easy to see that $V=\frac{\partial}{\partial\theta_{1}}+\frac{\partial}{\partial\theta_{2}}$ commutes with $\Delta_{sR}$, hence $\Delta_{sR}$ commutes with $\Delta_{S^{3}}$. Thus $\Delta_{sR}$ and $\Delta_{S^{3}}$ have a \emph{common} complete orthonormal set of eigenfunctions \cite{R:Dirac},\cite{R:vonNQM}; the eigenfunctions of $\Delta_{sR}$ are simply the spherical harmonics.

\sgap

Particularly noteworthy is $(x_{1}+iy_{1})^{k}=\sin^{k}\!\theta_{0}\,e^{ik\theta_{1}}$. It is a ``Gaussian beam'': a family of eigenfunctions of both $\Delta_{S^{3}}$ and $\Delta_{sR}$ that concentrates along a great circle. Zelditch singles out this example in the Riemannian setting (see, e.g., \cite{R:ZelditchParkCity}, p.76). It would be interesting to see if it is possible to construct, localized to each sR geodesic, a quasimode or Gaussian beam in the spirit of Ralston \cite{R:RalstonConstr}, \cite{R:RalstonApprox}.

\sgap

Taylor used the Peter-Weyl Theorem to find the eigenvalues of $\Delta_{sR}$ \cite{R:TaylorNonc}; Domokos generalized, using subelliptic Peter-Weyl and Plancherel Theorems on compact, connected, semisimple Lie groups \cite{R:DomokosPeterWeyl}. To summarize, the eigenvalues of $-\Delta_{S^{3}}$ are $m(m+2)$ for $m\in\{0,1,2,\ldots\}$, and the eigenvalues of $-\Delta_{sR}$ are (for the same $m$; the operators have the same complete orthonormal set of eigenfunctions):
$$4mj-4j^{2}+2m,\qquad j\in\{0,1,2,\ldots,m\}.$$ 
For reference, the eigenvalues of the $\lambda$-penalty Laplacian
$$-\Delta_{\lambda}=-\Delta_{sR}-\lambda^{-2}V^{2}$$ are
$$(1-\lambda^{-2})4j(m-j)+m(2+\lambda^{-2}m),$$
for $m\in\{0,1,2,\ldots\}$ and $j\in\{0,1,2,\ldots,m\}$.

\sgap

At this point we will not conjecture a general formula relating the sR length spectrum of a bracket-generating compact sR manifold (which for $S^{3}$ is $\{2\pi \sqrt{n};\,\,n\in\mathbb{N}\}$) to the set of eigenvalues of the subLaplacian counted with or without multiplicities (which for $S^{3}$ is $\{2m;\,\,m=0,1,2,\ldots\}$). 

\sgap

\end{document}